\documentclass[twoside]{amsart} 

\usepackage{amssymb, amsmath, amsfonts, eepic, graphicx}
\usepackage{pstricks,pst-poly,pst-node,pst-text,pst-plot,pstricks-add}
\usepackage{exscale}

\newcommand{\beq}{\begin{equation}}
\newcommand{\eeq}{\end{equation}}

\newtheorem{rem}{Remark}

\newtheorem{examp}{Example}

\begin{document}

\begin{center} {\LARGE  {\bf A Remark on Deriving Precise Upper Bounds of the Number of the RNA Secondary Structures  \vspace{.25cm} }  }  \\

{\Large    \emph{Alexander I. Kheyfits  \vspace{.25cm}   } }

akheyfits@gc.cuny.edu   \vspace{.25cm}

\end{center}

\noindent {\bf {Abstract}.} An elementary, at the undergraduate level, derivation is given of precise upper bounds of the number of the various RNA structures. The method works when the generating function has multiple singularities at the circle of convergence. Nor does the method require the Taylor coefficients to be positive. Examples from the current research literature are considered.
\vspace{.5cm}

\begin{center}
\today
\end{center}

\vspace{.5cm}

\footnoterule \vspace{.25cm}


\section{Introduction}

Since Waterman and Stein \cite{Wat, StWat} defined the RNA secondary structures in graph-theoretical terms, derivation of the upper bounds or precise asymptotic formulas for the various structures became an important problem; the number of relevant papers is growing, see, e.g., \cite{CKK, FuC, LoPoCl} and the references therein. Derivation of these estimates in the current literature is based on the deep result of complex analysis, that can be traced back to G. Darboux; we refer the reader to the very informative survey of various methods \cite{LoPoCl}.

During the workshop \emph{Teaching Discrete and Algebraic Mathematical Biology to Undergraduates} at the Mathematical Biosciences Institute at Ohio State University, Columbus, OH, 7/29/2013 - 8/02/2013, several speakers discussed these estimates. However, the Darboux theorem and its modern analogs are well beyond the current undergraduate curriculum.

The goal of this methodical note is to show that precise upper bounds for the number of the secondary structures in many cases can be derived quite elementary, within the power of an undergraduate student taking an Introductory Complex Analysis class, and without the use of CAS. Of course, the method is not limited to the secondary structures, in the last two examples we find the bounds for the number of RNA shapes \cite{LoPoCl}. The method is based on the well-known Cauchy-Hadamard formula for the radius of convergence of Taylor series, or even on its real-valued relative -- the root-test for convergence of the power series.

Moreover, unlike some other methods, it is irrelevant for the method, whether the generating function has more than one singularity on the circle of convergence. Nor does the method require the Taylor coefficients to be positive.

For the readers convenience, in the next section we remind what are the RNA primary and secondary structures. Then we describe our method, and in Section 3 we consider examples of its application.

\section{RNA Structures and the Cauchy-Hadamard Formula }

\subsection{RNA Primary Structures}

Living cells contain important (macro)molecules, called RiboNucleic Acids (RNA). These acids contain essential genetic information, for example, about  viruses, thus, it is important to know their structure. Biologists distinguish primary, secondary, and tertiary structures of RNA. Unlike the double helix of the DNA, each RNA is a  \emph{linearly ordered strand}, or just a \emph{string}, consisting of other molecules, called \emph{ribonucleotides}. This string is the backbone of any RNA. Traditionally, it is represented by a horizontal straight segment with nodes occupied by the nucleotides. If the RNA contains $n$ ribonucleotides, we select $n$ points of the segment and number them consecutively from left to the right by the natural numbers $1,2,\ldots ,n$; these dots represent the ribonucleotides in the RNA. For more information the reader can consult, for example, \cite{Re} or \cite{WaZ}.

When we start studying new objects, it is often necessary to know their quantity. In particular, it is important to know the number of the primary and secondary structures of RNA. It is not always possible to find a precise formula for the number of the secondary structures subject to various restrictions. And even if such a formula is derived, it can be very cumbersome, and therefore useless. That is why, a lot of work has been done to derive different \emph{asymptotic} formulas for the numbers of various secondary structures, see, e.g., \cite{Re, WaZ} and the references therein. A formula is called asymptotic, if it gives better and better \emph{relative approximation} of a quantity under consideration, when a certain parameter (e.g., the size of a system or time) is approaching a crucial threshold; for example, if time tends to infinity.

There are four different ribonucleotides, called \emph{adenine} (denoted hereafter $\mathcal{A}$), \emph{cytosine} ($\mathcal{C}$), \emph{guanine} ($\mathcal{G}$), and \emph{uracil} ($\mathcal{U}$). The linear ordering of these four nucleotides in either order, where each of them can repeat indefinitely, is called the \emph{primary structure} of the RNA. Thus, the primary structure of an RNA can be depicted by drawings like this,
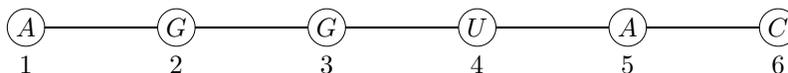
\begin{figure}[hbt]
\psset{xunit=1cm,yunit=1cm}
\def\xlim{6.5}
\def\ylim{.75}
\begin{pspicture*}(-\xlim,-\ylim)(\xlim,\ylim)
\psline(-4.75,0)(-3.25,0)      \psline(-2.75,0)(-1.25,0)   \psline(-.75,0)(.75,0)        \psline(1.25,0)(2.75,0)
\psline(4.75,0)(3.25,0)

\rput(-5,0){$\circle{.5}$}  \rput(-3,0){$\circle{.5}$}    \rput(-1,0){$\circle{.5}$}     \rput(1,0){$\circle{.5}$}   \rput(3,0){$\circle{.5}$} \rput(5,0){$\circle{.5}$}

\rput(-5,0){$A$}    \rput(-3,0){$G$}    \rput(-1,0){$G$}     \rput(1,0){$U$}  \rput(3,0){$A$}   \rput(5,0){$C$}

\rput(-5,-.5){$1$}    \rput(-3,-.5){$2$}    \rput(-1,-.5){$3$}     \rput(1,-.5){$4$}  \rput(3,-.5){$5$}   \rput(5,-.5){$6$}

\end{pspicture*}
\caption{A string (primary structure) with 6 nodes occupied by the nucleotides $\mathcal{A-G-G-U-A-C}$.}
\end{figure}

In mathematical parlance, the pictures like this are called \emph{graphs}. They are studied in the \emph{graph theory}. The graph in Fig. 1 is \emph{labeled} -- its vertices are labeled by the symbols of the nucleotides. The graph theory is a mathematical theory, and even though mathematics by itself cannot solve biological problems, it can give useful insights and help to solve biological problems \cite{Wig}.

\subsection{Counting the Primary Structures}

We begin by solving an easy problem of calculating the number of primary structures of RNA of some specified length, say $n$. We denote the number of different linear strings, containing $n$ nucleotides, without any restrictions on the neighboring ones, as $\widetilde{R}(n)$. Since every string starts with one of the four nucleotides, either $\mathcal{A}$, or $\mathcal{C}$, or $\mathcal{G}$, or $\mathcal{U}$, followed by a string of length $n-1$, we can immediately produce the basic equation
\[\widetilde{R}(n)=4\times \widetilde{R}(n-1).\]
Such equations are called \emph{recurrence relations} or \emph{difference equations}. In the same fashion,
\[\widetilde{R}(n-1)=4\times \widetilde{R}(n-2),\]
and we can \emph{iterate} this equation, getting the equation
\[\widetilde{R}(n)=4 \widetilde{R}(n-1)=4^2 \widetilde{R}(n-2)=4^3 \widetilde{R}(n-3)=\cdots =4^{n-1} \widetilde{R}(1).\]
Since we have an obvious initial condition $\widetilde{R}(1)=4$, the total number of primary structures without any restriction is
\[\widetilde{R}(n)=4^n.\]
We can notice that this is just the number of permutations (or arrangements) with repetitions of $n$ elements of four different kinds of elements\footnote{For all basic information from combinatorics and graph theory see for example \cite{Kh}.}.

Thus, the number of the RNA grows exponentially, as $4^n$. Therefore, the number of the secondary structures in the literature is usually compared with the exponential function $b^n$.

Now we take up the RNA strings with restrictions on the neighboring nucleotides. Of course, in real molecules there are always small deviations from the basic rules, that is, certain "forbidden" pairs can occur, even though with a small probability. We neglect these "outliers" and consider only RNA, where \emph{all the base pairs} are only of these three types, called Watson-Crick pairing:
\beq \mathcal{A-U;\; C-G;\; G-U}.\eeq

We compute the number of the primary structures satisfying these restrictions. The very first nucleotide can be either $\mathcal{A}$, or $\mathcal{C}$, or $\mathcal{G}$, or $\mathcal{U}$. Let us denote the number of strings of length $n$ and starting with $\mathcal{A}$, as $R_{\mathcal{A}}(n)$, and similarly, $R_{\mathcal{C}}(n)$, $R_{\mathcal{G}}(n)$, $R_{\mathcal{U}}(n)$. If the very first molecule is an $\mathcal{A}$, then the second nucleotide can be only $\mathcal{U}$. Since this $\mathcal{U}$ starts a string of $n-1$ bases,  thus
\[R_{\mathcal{A}}(n)=R_{\mathcal{U}}(n-1).\]
Similarly,
\[R_{\mathcal{C}}(n)=R_{\mathcal{G}}(n-1).\]
However, if the first nucleotide is a $\mathcal{G}$, then the second molecule must be either $\mathcal{C}$ or $\mathcal{U}$, thus
\[R_{\mathcal{G}}(n)=R_{\mathcal{C}}(n-1)+R_{\mathcal{U}}(n-1).\]
Similarly, if the first nucleotide is a $\mathcal{U}$, then the second molecule is either $\mathcal{A}$ or $\mathcal{G}$.

Since an RNA must start with some nucleotide, then
\[R(n)=R_{\mathcal{A}}(n)+R_{\mathcal{C}}(n)+R_{\mathcal{G}}(n)+R_{\mathcal{U}}(n).\]

Collecting all these equations together, we deduce the recurrence equation
\[R(n)=R_{\mathcal{U}}(n-1)+ R_{\mathcal{G}}(n-1)+ R_{\mathcal{C}}(n-1)+R_{\mathcal{U}}(n-1)+ R_{\mathcal{A}}(n-1)+R_{\mathcal{G}}(n-1)\]
\beq =R_{\mathcal{A}}(n-1)+R_{\mathcal{C}}(n-1)+2R_{\mathcal{G}}(n-1)+2R_{\mathcal{U}}(n-1). \eeq

To have the unique solution, we must supply the \emph{initial conditions}, which in this case are, obviously, \beq R_A(1)=R_C(1)= R_G(1)=R_U(1)=1;\; R(1)=4.\eeq

Using equations (2)-(3), we can easily compute the number $R(n)$ for any given $n$; for large $n$ we should probably use computers. For example, if $n=2$ then $n=14$ we get $R(2)=1+1+2+2=6$. Indeed, we can list these strands explicitly,
\[\mathcal{A-U};\; \mathcal{C-G};\; \mathcal{G-C};\; \mathcal{G-U},\; \mathcal{U-A};\; \mathcal{U-G}.\]

\subsection{Secondary Structures}

An RNA molecule is not rigid like a metal bar, it is flexible and can be conveniently thought of as a smooth flexible string, which can be crumpled, and then stretched again without any noticeable change.

Imagine now that we attached small pieces of velcro tape at some places of this string. If we now fold it over, then these pieces of velcro tape can hook one another, and we cannot easily stretch the tape in a linear structure as before. In real molecules instead of velcro tape there are certain pairs of nucleotides. If they happen to be close enough one to another, they are capable of forming chemical bonds. According to Watson-Crick, these are three \emph{base pairs} of the nucleotides, given by equation (1) in either order. For example, the string in Fig. 1 can fold over into the one in Fig. 2, where the new ties are pairs $\mathcal{C-G}$ and $\mathcal{G-C}$, but it cannot fold into the one in Fig. 3.
\begin{figure}[hbt]
\psset{xunit=1cm,yunit=1cm}
\def\xlim{6.5}
\def\ylim{2}
\begin{pspicture*}(-\xlim,-\ylim)(\xlim,\ylim)
\psline(-3.75,-1)(-2.25,-1)          \psline(-1.75,-1)(-.25,-1)      \psline(.25,-1)(1.75,-1)

\rput(-4,-1){$\circle{.5}$}    \rput(-2,-1){$\circle{.5}$}  \rput(0,-1){$\circle{.5}$}    \rput(2,-1){$\circle{.5}$}     \rput(2,1){$\circle{.5}$}   \rput(0,1){$\circle{.5}$}

\psline(2,-.75)(2,.75)        \psline(1.75,1)(.25,1)      \psline(.05,-.77)(.05,.73)    \psline(-.05,-.77)(-.05,.73)

\rput(-4,-1){$A$}    \rput(-2,-1){$G$}    \rput(0,-1){$G$}     \rput(2,-1){$U$}
\rput(2,1){$A$}   \rput(0,1){$C$}

\end{pspicture*}
\caption{A secondary structure built on the primary string $\mathcal{A-G-G-U-A-C}$.}
\end{figure}
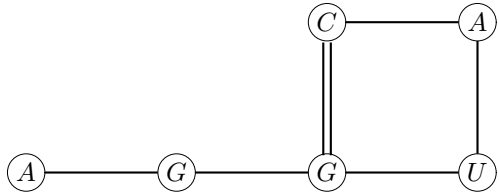

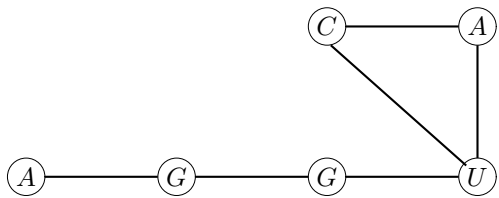
\begin{figure}[hbt]
\psset{xunit=1cm,yunit=1cm}
\def\xlim{6.5}
\def\ylim{2}
\begin{pspicture*}(-\xlim,-\ylim)(\xlim,\ylim)
\psline(-3.75,-1)(-2.25,-1)          \psline(-1.75,-1)(-.25,-1)      \psline(.25,-1)(1.75,-1)

\rput(-4,-1){$\circle{.5}$}    \rput(-2,-1){$\circle{.5}$}  \rput(0,-1){$\circle{.5}$}    \rput(2,-1){$\circle{.5}$}     \rput(2,1){$\circle{.5}$}   \rput(0,1){$\circle{.5}$}

\psline(2,-.75)(2,.75)        \psline(1.75,1)(.25,1)      \psline(.05,.75)(1.85,-.85)
\rput(-4,-1){$A$}    \rput(-2,-1){$G$}    \rput(0,-1){$G$}     \rput(2,-1){$U$}
\rput(2,1){$A$}   \rput(0,1){$C$}

\end{pspicture*}
\caption{This structure is forbidden by the Watson-Crick pairing.}
\end{figure}

This folding of RNA molecules in the plane is called the \emph{secondary structure} of the RNA molecule. Since different nucleotides can come close to each other, a primary structure can generate different secondary structures, which significantly complicates the analysis of the RNA. It is supposed that the secondary structure is a two-dimensional object, thus it can be drawn in the plane.

The \emph{secondary structure} of the RNA describes the ordering and location of these base pairs of the nucleotides. The secondary structure is responsible for many crucial biological phenomena, and the graph theory is helpful in discovering these structures.

The secondary structure, as defined by M. Waterman \cite{Wat}, ia also a graph. We consider hereafter only \emph{simple} graphs, i.e., graphs without loops or parallel edges. To a simple graph, there corresponds a square matrix of size  $n\times n$, such that its element $a{i,j}=1$ if and only if the vertices $v_i$ and $v_j$ are connected with an edge and is $0$ otherwise; this matrix is called the \emph{adjacency matrix} of the graph. The secondary structure is a horizontal backbone of length $n$, i.e., just a primary structure, enriched by several arcs in the upper half-plane, whose end-points are the nodes of the backbone. In notation and terminology we follow \cite{Re}. The arcs of the secondary structure are subject to certain restrictions. The arc with ends at the nodes $i$ and $j>i$ is denoted as $(i,j)$, the \emph{length} of the arc is $j-i\geq 1$. Not every family of arcs corresponds to a secondary structure.

A secondary structure is a  simple graph on the backbone of the length $n$, such that the adjacency matrix $A=(a_{i,j})$ possesses the following three properties.  \\

1) For the basic string to be a backbone, we require that $a_{1,2}=a_{2,3}=\cdots =a_{n-1,n}=1$.

2) Next, not counting the neighbors, every point can be adjacent to at most one other point of the backbone. In terms of the incidence matrix, this means that for any $i,\; 1\leq i\leq n,$ there exists at most one $j$ with $j\neq i\pm 1$, such that $a_{i,j}=1$.

3) Finally, it is assumed that if the vertices $a_i$ and $a_j,\; i<j,$ are adjacent and $i<k<j$, then the vertex $a_k$ \emph{cannot be adjacent} with any vertex to the left of $a_i$ or to the right of $a_j$. In terms of the adjacency matrix this means that if $a_{k,l}=1$ and $i < k < j $, then also $i < l < j $. \\

It follows that the arcs of a secondary structure \emph{do not intersect}, it is a \emph{non-crossing} structure.  \\

\section{Asymptotic Enumeration of the Secondary Structures. Examples}

\subsection{Convolution and Generating functions of Sequences}

Even if $n$ is about 10, the total listing of all the secondary structures is cumbersome, so that we want to estimate their quantity. A convenient device for  this is their \emph{generating function}, that is, the power series
\[S(x)=\sum_0^{\infty}S_n x^n,\]
where $S_n$ is the number of secondary structures with $n$ nods, $n=0,1,2,\ldots $, and $x$ is an indeterminate, real or complex. If the series is divergent, it can be considered as formal power series, or we can truncate the series and consider \emph{generating polynomials}, as discussed, e.g., in \cite{Kh}. However, in the following examples all the series have positive radii of convergence.

Given two sequences, $a=\{a_0,a_1,\ldots \}$ and $b=\{b_0,b_1,\ldots \}$ with generating functions $f_a(x)$ and $f_b(x)$ respectively, the generating function of their linear combination $\alpha a+\beta b$ with any real or complex coefficients $\alpha$ and $\beta$ is the linear combination
\[f_{\alpha a+\beta b}(x)=\alpha f_a(x)+\beta f_b(x),\]
thus, the correspondence between the sequences and their generating functions is a linear transformation. There is also an operation with sequences, which corresponds to the multiplication of their generating functions. Indeed, let $f_a(x)$ and $f_b(x)$ be two absolutely convergent power series. Multiplying them termwise and combining like terms, we get the series
\[f_a(x)\times f_b(x)=\sum_0^{\infty} c_n x^n, \]
where
\[c_0=a_0\cdot b_0,\; c_1=a_0\cdot b_1+a_1\cdot b_0,\ldots , c_n=a_0\cdot b_n+ a_1\cdot b_{n-1}+a_2\cdot b_{n-2}+\cdots +a_n\cdot b_0,\ldots .\]

Thus defined sequence
\[c=\{c_n\}_{n=0}^{\infty}\]
is called the \emph{convolution} or the \emph{Cauchy product} of the sequences $a$ and $b$, and corresponds to the multiplication of the generating functions.  \\

Let $S^{\lambda}(x)$ be the generating function of the secondary structures with the arc length at least $\lambda \geq 2$. If $\lambda=2$, the recurrent relation
\beq S^{\lambda}_n=S^{\lambda}_{n-1}+\sum^{n-1-\lambda}_{j=0} S^{\lambda}_{n-2-\lambda} S^{\lambda}_j \eeq
was derived by Waterman \cite{Wat}. To solve it, we must supply $\lambda+1$ initial conditions. We assume
\[S^{\lambda}_0=S^{\lambda}_1=\cdots =S^{\lambda}_{\lambda}=1.\]

\subsection{Computing the Generating Functions for the Secondary Structures}

The crucial observation is that recurrent relation (4) contains a sum of pairwise products quite similar to the equation for the convolution. First, consider the case $\lambda =2$ and to simplify writing, set $S(x)=S^{\lambda}(x)$. Multiplying (4) by $x^n$, after some simple algebra one derives the following quadratic equation\footnote{However, the method has the broader scope. In Example 3 below, the method is applied to the generating function satisfying a cubic equation.} for the generating function,
\beq x^2 S^2(x)+(x-1-x^2)S(x)+1=0. \eeq
Solving it by the quadratic formula, we find
\[S(x)=\frac{1}{2x^2}\left(x^2-x+1\pm \sqrt{1-2x-x^2-2x^3+x^4} \right). \]

From (5) we see that $S(0)=1$, therefore, we have to choose the sign $"-"$ above, and finally get
\[S(x)=\frac{1}{2x^2}\left(x^2-x+1-\sqrt{1-2x-x^2-2x^3+x^4} \right). \]

\subsection{Estimations of the Number of the Secondary Structures}

We need the following facts from introductory Complex Analysis course. If a power series $f(z)=\sum_{n=0}^{\infty}f_nz^n$ has a positive radius of convergence $R$, then its sum $f(z)$ is an analytic function in the open disc $|z|<R$, $R$ is the distance from the $z=0$ to the nearest singular point, which must be located on the boundary $|z|=R$, and the Cauchy-Hadamard formula
\[\frac{1}{R}=\limsup_{n\rightarrow \infty}\sqrt[n]{|f_n|} \]
is valid. We can even (with some reservations) refer, instead of the Cauchy-Hadamard formula, to the root test, which is more-or-less real valued version of the latter.

Therefore, $\sqrt[n]{|f_n|}\leq 1/R$, and we have the required upper bound for the coefficients $f_n$,
\beq |f_n|\leq \left(\frac{1}{R}\right)^n,\; n=1,2,3,\ldots . \eeq

By the definition of the upper limit, there is a subsequence $f_{n_k}\rightarrow \infty$ as $k\rightarrow \infty$, such that \[\frac{1}{R}=\limsup_{n\rightarrow \infty}\sqrt[n]{|f_n|}. \]
This shows that the upper bound (6) is precise, it cannot be made smaller. Moreover, the method works if there are several singularities on the circle of convergence. Nor does the method require the Taylor coefficients to be positive. However, preserving this elementary level, we can derive only the precise exponential upper bound for the quantities at question, and cannot find the principal term of the asymptotic formulas explicitly, as more sophisticated methods do. To demonstrate the method, we apply it to several known examples from the existing literature.

\begin{examp} Consider the generating function $S(x)$ above. The origin $x=0$ is a \emph{removable singularity}, as can be straightforwardly seen by rationalizing the numerator,
\[S(x)=\frac{2}{x^2-x+1+\sqrt{1-2x-x^2-2x^3+x^4}}. \]
The radicand $1-2x-x^2-2x^3+x^4$ is a symmetric polynomial of forth degree. The corresponding equation $x^4-2x^3-x^2-2x+1=0$ can be explicitly solved by dividing over $x^4$ and substituting $t=x+1/x$; it has two real roots, $z_1=\frac{3-\sqrt{5}}{2}$, $z_2=\frac{3+\sqrt{5}}{2}$, and two complex conjugate roots, $z_{3,4}=\frac{-1\pm \imath \sqrt{3}}{2}$; $|z_{3,4}|=1$. The root, closest to the origin, is $z_1$, thus $S(z)$ is analytic in the open disc $|z|<R=|z_1|<1$. Finally, from (6) we get the estimate
\[|S_n|\leq \left(\frac{2}{3- \sqrt{5}}\right)^n =\left(\frac{3+ \sqrt{5}}{2} \right)^n,\]
which is, up to the pre-exponential factor, the estimate derived in \cite[p, 65]{Re}, \cite{LoPoCl}, and some other places  by making use of certain more sophisticated tools of complex analysis.  \end{examp}

\begin{rem} It must be noted that the function $S(z)$ is a many-valued analytic function of the complex argument $z$ with the branching points at the roots $z_1-z_4$ of the radicand and at infinity. To select a single-valued branch of $S$, we can cut the plane, for example, along the rays going from the points $z_1, z_3, z_4$ to infinity. The generating function $S(x)$ is given by the principal branch of $S(z)$ corresponding to the choice $\sqrt{1}=+1$. In the disk $|z|<|z_1|$, this generating function $S(x)$ is represented by the Taylor series with positive coefficients $S_n$. The same remark is valid for all the examples that follow. \end{rem}

\begin{examp}  In general case, that is, for the secondary structures with any arc length $\lambda \geq 2$, the difference equation is (\cite{Re})
\[S^{\lambda}(n)=S^{\lambda}(n-1)+\sum^{n-1-\lambda}_{j=0} S^{\lambda}(n-2-j) S^{\lambda}(j),\]
leading again to a quadratic equation for the generating function,
\[x^2 \left(S^{\lambda}(x)\right)^2-(1-x+x^2+\cdots +x^{\lambda})S^{\lambda}(x)+1=0, \]
which again can be explicitly solved for $S^{\lambda}$ by the quadratic formula. For example, if $\lambda=3$, the radicand is the polynomial $1-2x-x^2-x^4+2x^5+x^6$, which can be factored,
\[1-2x-x^2-x^4+2x^5+x^6=(1-2x-x^2)(1-x^4),\]
with the smallest root $\sqrt{2}-1$. Repeating the same reasoning as above, we get the estimate
\[S^3(n) \leq \left(\sqrt{2}+1\right)^n. \]
If $\lambda =4$, the radicand is a polynomial of $8^{th}$ degree, whose roots are to be evaluated numerically; the smallest one is $~0.436911$, leading to the estimate
\[S^4(n) \leq 2.28879^n. \] \end{examp}

\begin{examp}  Many authors consider secondary structures subject to various restrictions, see for example, \cite{WaZ, WZW, Kra} and the references therein. The presented approach works even if the equation for the generating function is of degree higher than $2$, as in the next example. We consider \emph{saturated secondary structures} studied in \cite{Kra}. In this case the generating function $S=S(z)=\sum_{n\geq 0}S_n z^n$ satisfies the system of two nonlinear equations (see, e.g., \cite{Kra} and the references therein)
\[S(z)=z+z^2+zT(z)+z^2T+z^2S+z^2S^2, \]
\[T(z)=z^2S+z^2TS.\]
Eliminating $T$, one derives the cubic equation for $S$,
\[z^4S^3+z^2(z^2-2)S^2+(1-z^2)S-z(1+z)=0. \]
We solve the cubic equation for $S$ by the classical Cardano formula. When the parameter $z$ is within the range of interest, the equation has one real root,
\[S_r(z)=\frac{2-z^2}{3z^3}+\frac{1}{3\sqrt[3]{2}z^3}\left\{A^{1/3}+
\frac{\sqrt[3]{2}(z^4-z^2+1)}{3A^{1/3}}\right\},\]
where
\[A=-2z^6+30z^4+27z^3+3z^2-2 \vspace{.3cm} \]
\[+(3z)^{3/2}\sqrt{-4z^7-4z^6+32z^5+60z^4+35z^3+6z^2-5z-4}.\]
The equation $A=0$ simplifies to $(z^4-z^2+1)^3=0$ with all the roots on the unit circle. The radicand in the $A$ has the smallest root at $z_0\approx0.424687310$, which is the radius of convergence, $R$, of $S(z)$. It should be mentioned that we chose the branch of the radical which is positive in a right neighborhood of $z_0$. Hence, we get the bound
\[S_n\leq \cdot (1/R)^n =const \cdot 2.354673^n \]
in complete agreement with \cite{Kra}.  \end{examp}

\begin{examp}  In the case of \emph{canonical secondary structures} \cite{Kra}, the generating function $S$ also satisfies a system of two nonlinear equations,
\[S(z) = z + zS(z) + z^2Q(z) + z^2S(z)Q(z),\]
\[Q(z) = z^3 + z^2Q(z) + z^4S(z)Q(z) + z^3S(z).\]
Eliminating $Q$, we get the quadratic equation for $S$ with the discriminant
\[\Delta(z)=z^{10}-4z^9-2z^8+6z^7+3z^6-8z^5-z^4+4z^3-z^2-2z+1=0.\]
Its smallest in absolute value root is $R\approx 0.5081360362$, which is the closest to the origin singular point of $S(z)$, and the upper bound for the number of the saturated secondary structures on $n$ nucleotides is its reciprocal $1/R~1.967977$, again in perfect agreement with \cite{Kra}.   \end{examp}

\begin{examp} Next, we consider an example from \cite[Theor. 2, p. 352]{FuC}, where the authors derived the generating function for a certain class of secondary structures,
\[S_0(z)=\frac{1}{2z^2(1+z)^2}\left(1-z-z^2(1+z)^2-\sqrt{P(z)} \right) \vspace{.2cm}\]
\[P(z)=(z^4+2z^3+z^2+z-1)^2-4z^3(1+z)^2.\]
Here $z=0$ is a removal singularity and $z=-1$ is beyond the disc of convergence, since the closest to the origin singular point is the branching point of the radical at the smallest root of the polynomial $P$. This root is $z_0\approx 0.32471796$, thus the radius of convergence is $R\approx 3.0795963$ and the upper bound of the number of the secondary structures at question is $\cdot 3.0795963$, again in agreement with \cite{FuC}.         \end{examp}

In the two examples to follow, we bound the number of various RNA shapes. The equations for the generating functions were derived in the informative paper by Lorenz, Ponty, and Clote \cite{LoPoCl}, we only suggest an elementary derivation of the main exponential term of the asymptotics.

\begin{examp} In the case of $\pi-$shapes, the generating function $S(z)$ satisfies the simple quadratic equation
\[S=z^2S^2+z^2S+z^2,\]
therefore, since $S(0)=0$,
\[S(z)=\frac{1}{2z^2}\left(1-z^2-\sqrt{1-2z^2-3z^4} \right).\]
The radicand is a biquadratic equation with roots $\pm \imath$ and $\pm 1/\sqrt{3}$. Thus, the distance from the boundary to the nearest singularity is $1/\sqrt{3}$, and the bound we sought for, is $s_n\leq (\sqrt{3})^n$.

It should be emphasized that it is irrelevant for this method whether there are multiple singularities on the circle of convergence.
\end{examp}
\begin{examp} In the last example we estimate the number of $\pi-$shapes compatible with the RNA sequences of length $n$. The equation for the generating function, also taken from \cite{LoPoCl}, is
\[z^2(1-z)^2S^2+(z-1+z^5-z^6)S+z^5=0.\]
Solving it by the quadratic formula and canceling by $1-z$, we get
\[S(z)=\frac{1}{2z^2(1-z)}\left(1-z^5-\sqrt{z^10-4z^7-2z^5+1} \right).\]
The singularity, closest to the origin, which gives the radius of analyticity of the generating function, is $R \approx ).756328$, hence, the required upper bound is its reciprocal $1/R \approx 1.32218$.
\end{examp}

\subsection*{\textbf{Acknowledgement}} The author is thankful to the MBI for inviting him to the workshop, and to Professors Raina Robeva and Peter Clote for useful discussions. Numerical computations in the examples above were made by making use of the Wolfram Mathematica free online root finder. The author wants to acknowledge the excellent performance of this widget. \\

\end{document}